\documentclass[a4paper,leqno]{amsart}
\usepackage{amsmath, amssymb, amsfonts, euscript} 

\DeclareMathOperator{\st}{|} 
\DeclareMathOperator{\otensprod}{\otimes} 
\DeclareMathOperator{\ctensprod}{\overline{\otimes}}
\DeclareMathOperator{\spec}{\sigma} 
\DeclareMathOperator{\Lip}{Lip}

\renewcommand{\geq}{\geqslant}
\renewcommand{\leq}{\leqslant}

\newcommand{\R}{\mathbb{R}}
\newcommand{\N}{\mathbb{N}}

\providecommand{\abs}[1]{\lvert #1 \rvert}
\providecommand{\scal}[2]{\langle #1,#2 \rangle}
\providecommand{\norm}[1]{\lVert #1 \rVert}
\providecommand{\set}[1]{\left\{ #1 \right\}}

\newcommand{\cd}[1]{\mathop{\bf #1}\nolimits}

\newcommand{\res}[1]{_{\rceil #1}}
\newcommand{\adh}[1]{\, \overline{#1}}

\def\dis{\displaystyle}
\def\eps{\varepsilon}
\def\Tbar{\bar{T}} 
\def\Tdbar{\underline{T\!\!}\,\,} 
\def\dbar{\bar{d}}
\def\ddbar{\underline{d}}
\def\Mt{\tilde{M}}
\def\nt{\tilde{n}}
\def\Ct{\tilde{C}}

\def\loc{_{\text{loc}}}

\def\V{\mathcal{V}}
\def\E{\mathcal{E}}
\def\O{\mathcal{O}}
\def\U{\mathcal{U}}
\def\A{\mathcal{A}}
\def\B{\mathcal{B}}
\def\CC{\mathcal{C}}

\def\Deltat{\tilde{\Delta}}

\providecommand{\scd}[2]{\scal{#1}{#2}}

\def\dxi{\dot{\xi}}
\def\dz{\dot{z}}

\def\dzeta{\dot{\zeta}}

\def\L{\mathcal{L}} 
\def\D{D} 

 \theoremstyle{plain}
 \newtheorem{theorem}{Theorem}[section]
 \newtheorem{theorem*}{Theorem}
 \newtheorem{lemma}[theorem]{Lemma}

 \theoremstyle{definition}
 \newtheorem{definition}[theorem]{{\sc Definition}}

 \theoremstyle{remark}
 \newtheorem{remark}[theorem]{{\sc Remark}}

\begin{document}


\title[null-controllability of the heat equation]{
On the null-controllability of the heat equation in unbounded domains
}
\author[L. Miller]{Luc Miller}
\address{
{\'E}quipe Modal'X, EA 3454 \\
Universit{\'e} Paris X, B{\^a}t.~G,
200 Av.~de la R{\'e}publique,
92001 Nanterre, France.}
\address{Centre de Math{\'e}matiques Laurent Schwartz, UMR CNRS 7640 \\ 
{\'E}cole Polytechnique, 91128
Palaiseau, France.}
\email{miller@math.polytechnique.fr}

\date{April 23, 2004. Accepted by the Bulletin des Sciences Math{\'e}matiques April 30, 2004.}

\begin{abstract}
We make two remarks about the null-controllability 
of the heat equation with Dirichlet condition in unbounded domains.
Firstly, 
we give a geometric necessary condition 
(for interior null-controllability in the Euclidean setting)
which implies that
one can not go infinitely far away from the control region
without tending to the boundary (if any),
but also applies when the distance to the control region is bounded.
The proof builds on heat kernel estimates.
Secondly, 
we describe a class of 
null-controllable heat equations on unbounded product domains.
Elementary examples include 
an infinite strip in the plane controlled from one boundary 
and an infinite rod controlled from an internal infinite rod.
The proof combines earlier results on compact manifolds
with a new lemma saying that 
the null-controllability of an abstract control system 
and its null-controllability cost 
are not changed by taking its tensor product 
with a system generated by a non-positive self-adjoint operator.
\end{abstract}






\maketitle

\tableofcontents 


\section{Introduction.}
\label{sec:intro}

\subsection{The problem.}
\label{sec:pb}

Let $M$ be a smooth connected complete 
$n$-dimensional Riemannian manifold 
with boundary $\partial M$. 
When $\partial M\neq\emptyset$, $M$ denotes the interior 
and $\adh{M}=M\cup\partial{M}$.
Let $\Delta$ denote the (negative) Laplacian on $M$.

Consider a positive control time $T$ 
and a non-empty open control region $\Gamma$ of $\partial M$.
Let ${\bf 1}_{]0,T[\times \Gamma}$ 
denote the characteristic function of the space-time 
control region $]0,T[\times \Omega$.
The heat equation on $M$ is said to be 
{\em null-controllable}
in time $T$ by boundary controls on $\Gamma$ 
if for all  
$\phi_{0}\in L^{2}(M)$
there is a control function 
$u\in L^{2}\loc(\mathbb{R};L^{2}(\partial M))$
such that the solution $\phi\in C^{0}([0,\infty),L^{2}(M))$
of the mixed Dirichlet-Cauchy problem:
\begin{gather} 
\label{eqHeat}
\partial_{t}\phi - \Delta \phi=0
\quad {\rm in}\ ]0,T[\times M, \quad 
\phi=\cd{1}_{]0,T[\times \Gamma} u   \quad {\rm on}\ ]0,T[\times\partial M,
\end{gather}
with Cauchy data
$\phi=\phi_{0}$ at $t=0$,
satisfies 
$\phi=0$ at $t=T$.
The {\em null-controllability cost} is the best constant, 
denoted $C_{T,\Gamma}$, in the estimate:
$$
\|u\|_{L^{2}(]0,T[\times \Gamma)}\leq C_{T,\Gamma}\|\phi_{0}\|_{L^{2}(M)}
$$
for all initial data $\phi_{0}$ and control $u$ 
solving the null-controllability problem described above.
The analogous interior null-controllability problem
from a non-empty open subset $\Omega$ of $\adh{M}$
is also considered:
\begin{gather}
\label{eqHeatI}
\begin{split}
&\partial_{t}\phi-\Delta\phi =\cd{1}_{]0,T[\times \Omega} u
\text{ on }\ \R_{t}\times M, 
\quad 
\phi=0\ \text{ on }\ \R_{t}\times \partial M, \\
&\phi(0)=\phi_{0}\in L^{2}(M),\ 
u\in L^{2}\loc(\R;L^{2}(M)) 
. \end{split}
\end{gather}

When $M$ is compact (for instance a bounded domain of the Euclidean space), 
Lebeau and Robbiano have proved 
(in \cite{LR95} using local Carleman estimates) 
that, for all $T$ and $\Gamma$ there is a continuous linear operator 
$S:L^{2}(M)\to C^{\infty}_{0}(\mathbb{R}\times \partial M)$
such that $u=S\phi_{0}$ yields 
the null-controllability of the heat equation (\ref{eqHeat}) on $M$ 
in time $T$ by boundary controls on $\Gamma$.
They have also proved the analogous result for (\ref{eqHeatI})
which implies that 
interior null-controllability holds for arbitrary $T$ and $\Omega$.
(We refer to \cite{FI96} for a proof of null-controllability  
for more general parabolic problems using global Carleman estimates.)

The null-controllability of the heat equation 
when $M$ is an unbounded domain of the Euclidean space 
is an open problem which Micu and Zuazua have recently underscored
in \cite{MZ03}.
On the one hand, it is only known to hold when 
$M\setminus\Omega$ is bounded (cf.~\cite{CMZ01}).
On the other hand, 
its failure can be much more drastic than in the bounded case 
(when $M$ is the half space and $\Gamma=\partial M$,
it is proved in \cite{MZ01a,MZ01b}
that initial data whith Fourier coefficients 
that grow less than any exponential 
are not null-controllable in any time,
whereas there are initial data with exponentially growing Fourier coefficients
that are null-controllable).

The geometric aspect of the open problem in \cite{MZ03} is addressed here
with examples of null-controllability with unbounded uncontrolled region,
and lack thereof 
including when the distance to the controlled region is finite
(cf.~th.\ref{th:rod}.iii).
The geometric necessary condition in th.\ref{th:GNC}
grasps at some notion of ``controlling capacity'' of a subset  
that would yield a necessary and sufficient condition
for interior null-controllability.

\subsection{Elementary examples.}
\label{sec:elem}

Before stating the results in full generality, 
we give elementary examples.

The simplest (bounded) case to study is 
when $M$ is a segment and $\Gamma$ is one of the end points.
It is well-known that this problem 
reduces by spectral analysis to classical results 
on non-harmonic Fourier series.
For further reference, we introduce 
the optimal fast control cost rate for this problem:
\begin{definition}
\label{def:alpha}
The rate $\alpha_{*}$ is the smallest positive constant such that 
for all $\alpha > \alpha_{*}$
there exists $\gamma>0$ such that,
for all $L>0$ and $T\in \left]0,\inf(\pi,L)^{2}\right]$,
the null-controllability cost $C_{L,T}$ 
of the heat equation (\ref{eqHeat}) on the Euclidean interval $M=]0,L[$ 
(i.e. $\Delta=\partial_{x}^{2}$) from $\Gamma=\{0\}$
satisfies:
$C_{L,T}\leq \gamma \exp(\alpha L^{2}/T)$.
\end{definition}
Computing $\alpha_{*}$ is an interesting open problem.
As proved in~\cite{LMheatcost}, 
\begin{theorem}
\label{th:1d}
The rate $\alpha_{*}$ defined above satisfies: 
$1/4\leq \alpha_{*}\leq 2\left(36/37\right)^{2}< 2$.
\end{theorem}

The simplest unbounded case where null-controllability holds 
is probably the following, 
which extends to an infinite strip 
the null-controllability 
from one side of a rectangle proved in \cite{Fat75}.

\begin{theorem}
\label{th:strip}
The heat equation (\ref{eqHeat})
on the infinite strip $M=]0,L[\times \R$ of the Euclidean plane
(i.e. $\Delta=\partial_{x}^{2}+\partial_{y}^{2}$)
is null-controllable 
from one side $\Gamma=\{ (x,y) | x=0, y\in \R \}$
in any time $T>0$.
Moreover, the corresponding null-controllability cost satisfies 
(with $\alpha_{*}$ as in th.\ref{th:1d}):
$\displaystyle
\limsup_{ T\to 0}  
T \ln C_{\Gamma,T} \leq \alpha_{*} L^{2}
$.
\end{theorem}

Here is an example in the usual three dimensional space 
which illustrates interior null-controllability and lack thereof.

\begin{theorem}
\label{th:rod}
Consider the heat equation (\ref{eqHeatI})
on the infinite rod $M=S\times \R$ in the Euclidean space
(i.e. $\Delta=\partial_{x}^{2}+\partial_{y}^{2}+\partial_{z}^{2}$)
where the section $S$ is any smooth connected bounded open set of the plane.

i) It is null-controllable in any time $T>0$
from any interior infinite rod $\Omega=\omega\times \R$ 
where the section $\omega$ is an open non empty subset of $\adh{S}$.
Moreover, if $\omega$ contains a neighborhood of the boundary of $S$
and $S\setminus\omega$ does not contain any segment of length $L$, 
then the corresponding null-controllability cost satisfies 
(with $\alpha_{*}$ as in th.\ref{th:1d}):
$\displaystyle
\limsup_{ T\to 0}  
T \ln C_{\Omega,T} \leq \alpha_{*} L^{2}
$.

ii) It is not null-controllable in any time $T>0$
from any interior region $\Omega$ 
of finite Lebesgue measure such that 
$M\setminus\Omega$ contains slabs $S\times [z_{1},z_{2}]$
of arbitrarily large thickness $\abs{z_{2}-z_{1}}$. 

iii) It is not null-controllable in any time $T>0$
from the cylindrical interior region 
$\Omega=\set{(x,y,z)\in M \st x^{2}+y^{2}< R(z)^{2}}$ 
if $(0,0)\in S$ 
and the lower semi-continuous 
function $R:\R\to [0,\infty)$ 
tends to zero at infinity.
\end{theorem}


\subsection{Main results.}
\label{sec:results}

A large class of null-controllable heat equations on unbounded domains 
is generated by the two following theorems 
concerning respectively boundary and interior controllability.
In both theorems, $\Mt$ denotes 
another smooth complete $\nt$-dimensional Riemannian manifold
and $\Deltat$ denotes the corresponding Laplacian.

\begin{theorem}
\label{th:prodgamma}
Let $\gamma$ denote the subset $\Gamma\times \Mt$ 
of $\partial(M\times\Mt)$.
If the heat equation (\ref{eqHeat}) is null-controllable
at cost $C_{T,\Gamma}$ then the heat equation: 
\begin{gather*}
\begin{split}
&\partial_{t}\phi-(\Delta+\Deltat)\phi =0\ 
\text{ on }\ \R_{t}\times M\times\Mt, 
\quad 
\phi={\bf 1}_{\gamma} g\ \text{ on }\ \R_{t}\times \partial (M\times\Mt), \\
&\phi(0)=\phi_{0}\in L^{2}(M\times\Mt),\ 
g\in L^{2}\loc(\R;L^{2}(\partial (M\times\Mt))) 
, \end{split}
\end{gather*}
is exactly controllable in any time $T$ 
at a cost $\Ct_{T,\gamma}$ which is not greater than $C_{T,\Gamma}$.
\end{theorem}

\begin{theorem}
\label{th:prodomega}
Let $\omega$ denote the subset $\Omega\times \Mt$ 
of $M\times\Mt$.
If the heat equation (\ref{eqHeatI}) is null-controllable
at cost $C_{T,\Omega}$ then the heat equation: 
\begin{gather*}
\begin{split}
&\partial_{t}\phi-(\Delta+\Deltat)\phi ={\bf 1}_{\omega} g
\text{ on }\ \R_{t}\times M\times\Mt, 
\quad 
\phi=0\ \text{ on }\ \R_{t}\times \partial (M\times\Mt), \\
&\phi(0)=\phi_{0}\in L^{2}(M\times\Mt),\ 
g\in L^{2}\loc(\R;L^{2}(M\times\Mt)) 
, \end{split}
\end{gather*}
is exactly controllable in any time $T$ 
at a cost $\Ct_{T,\omega}$ which is not greater than $C_{T,\Omega}$.
\end{theorem}

\begin{remark}
Th.\ref{th:rod} i) is a particular case of th.\ref{th:prodomega}
with $M=S$, $\Mt=\R$, inverted $\Omega$ and $\omega$,
and the cost estimate results from the cost estimate on $M$
proved in \cite{LMheatcost}.
Th.\ref{th:prodgamma} and th.\ref{th:prodomega}
apply, for instance, 
to any open subset $\Mt$ of the Euclidean space $\R^{\nt}$. 
Thanks to the results of \cite{LR95} 
already mentioned in section~\ref{sec:pb},
the conclusions of these theorems 
hold for arbitrary control regions of a compact $M$. 
Then they can be applied recursively, 
taking the resulting null-controllable product manifold 
as the new $M$ (the theorems are still valid if $M$ has corners).
\end{remark}

\begin{remark}
The case when $M$ is a bounded Euclidean set 
and $\Mt=(0,\eps)$ with Neumann boundary conditions at both ends
has been considered in \cite{TZ00}
with an extra time-dependent potential. 
When $\eps\to 0$, using global Carleman estimates, 
it is proved that the cost is uniform
(as in th.\ref{th:prodomega})
and depends on the uniform norm of the potential.
Moreover, the limit of the control functions is a control function 
for the limit problem.
\end{remark}

\begin{remark}
The type of boundary conditions are irrelevant 
to the proof of th.\ref{th:prodgamma} and th.\ref{th:prodomega}.
These theorems 
can be combined with 
th.6.2 in  \cite{LMcoscost}
and th.2.3 in \cite{LMheatcost} respectively 
to obtain 
bounds on the fast null-controllability cost:
$$
\limsup_{ T\to 0}  
T \ln \Ct_{\gamma,T} \leq \alpha_{*} L_{\Gamma}^{2}
\quad \text{ and } \quad
\limsup_{ T\to 0}  
T \ln \Ct_{\omega,T} \leq \alpha_{*} L_{\Omega}^{2}
$$
for any  
$L_{\Gamma}$ and $L_{\Omega}$ 
such that every generalized geodesic of length greater than $L_{\Gamma}$ 
passes through $\Gamma$ at a non-diffractive point,
and every generalized geodesic of length greater than $L_{\Omega}$
passes through $\Omega$.
We refer readers interested by these bounds 
to \cite{LMheatcost, LMcoscost} 
where more is said about generalized geodesics
and the extra geometric assumptions needed to use them. 
\end{remark}

The last result states a geometric condition 
which is necessary for the interior null-controllability of the heat equation 
on an unbounded domain of the Euclidean space.
This condition involves the following ``distances''.

\begin{definition}
\label{def:dist}
In $\R^{n}$, the Euclidean distance of points from the origin 
and the Lebesgue measure of sets are both denoted by $\abs{\cdot}$.
Let $M$ be a non-empty open subset of $\R^{n}$.
Let $d:\adh{M}^{2}\to \mathbb{R}_{+}$ denote the distance function on $M$,
i.e. the infimum of lengths of arcs in $M$ with end points $x$ and $y$
(n.b., in terms of Lipschitz potentials:
$d(x,y)=\sup_{\psi\in \Lip(\adh{M}), \norm{\nabla \psi}_{L^{\infty}}\leq 1}
\abs{\psi(x)-\psi(y)}$). 
The distance of $y\in M$ from the boundary of $M$
is $d_{\partial}(y)=\inf_{x\in \R^{n}\setminus M}\abs{x-y}$.
The distance of $y\in \adh{M}$ from $\Omega\subset M$ 
is $d(y,\Omega)=\inf_{x\in \Omega}d(x,y)$.
We define the 
{\em averaged distance} $\dbar_{T}(y,\Omega)$
of $y$ to $\Omega$ 
with Gaussian weight of variance $T$ by
\begin{gather*} 
\dbar_{T}(y,\Omega)^{2}=-2T 
\log\left( \int_{\Omega} \exp\left( -\frac{d(y,x)^{2}}{2T} \right) dx \right)
\geq d(y,\Omega)^{2}-2T\log\abs{\Omega} 
\ . 
\end{gather*}
Technically, we shall use the following 
{\em bounded distance} of $y$ to $\partial M$:
\begin{gather*} 
\ddbar_{T}(y,\partial M)=\min\set{d_{\partial}(y),T\pi^{2}n/4} \ .
\end{gather*}
\end{definition}

\begin{theorem}
\label{th:GNC}
Let $M$ be a connected open subset of $\R^{n}$
and let $\Omega$ be an open subset of $M$.
If there are a sequence $\set{y_{k}}_{k\in\N}$ of points in $M$,
a time $\Tbar>0$ and a constant $\kappa>1$ such that 
\begin{gather}
\label{eqgeomcond}
\dbar_{\Tbar}(y_{k}, \Omega)^{2}-
\kappa 
\frac{\pi^{2}n^{2}}{4} 
\left(\frac{\Tbar}{\ddbar_{\Tbar}(y_{k},\partial M) }\right)^{2}
\to +\infty\, , \text{ as } k\to +\infty 
\ ,
\end{gather}
then the heat equation (\ref{eqHeatI}) 
is not null-controllable in any time $T<\Tbar$.
In particular, when 
$\Omega$ has finite Lebesgue measure, if there is 
a sequence $\set{y_{k}}_{k\in\N}$ 
such that $\inf_{k}d_{\partial}(y_{k})>0$ 
and $\lim_{k}d(y_{k}, \Omega)=\infty$,
then the heat equation 
(\ref{eqHeatI}) 
is not null-controllable in any time $T$.
\end{theorem}

\begin{remark} 
The simple condition in the second part of th.\ref{th:GNC}
is enough to prove th.\ref{th:rod} ii)
(consider the points $(0,0,(z_{2}-z_{1})/2)$ 
of a sequence of slabs $S\times [z_{1},z_{2}]$ in $M\setminus\Omega$
with thickness $\abs{z_{2}-z_{1}}$ tending to infinity).
Th.\ref{th:rod} iii) 
illustrates that it may fail 
although the finer condition (\ref{eqgeomcond}) holds.
The second term in the geometric condition (\ref{eqgeomcond}) 
allows $\set{y_{k}}_{k\in\N}$ to tend to the boundary of $M$.
To illustrate its usefulness, we give yet another example 
in rk.\ref{rem:shrinkrod}. 
\end{remark}

\begin{remark} \label{rem:NCmanifolds}
The proof of th.\ref{th:GNC} in sect.\ref{sec:proof} 
builds on heat kernel estimates.
Generalizations to some non-compact manifolds 
can obviously be obtained 
using the heat kernel estimates available in the literature
(cf.~\cite{Zha03} and ref.~therein).
We consider null-controllability on non-compact manifolds
in a forthcoming paper.
\end{remark}


\section{An abstract lemma on tensor products}
\label{sec:tensor}

In this section, we prove that 
the cost of null-controllability of an abstract control system 
is not changed by taking its tensor product 
with an uncontrolled system generated 
by a non-positive self-adjoint operator.

\subsection{Abstract setting}
\label{sec:abstract}

We first recall 
the general setting for control systems:
admissibility, observability and controllability
notions and their duality (cf.~\cite{DR77} and~\cite{Wei89}).

Let $Z$ and $\V$ be Hilbert spaces.
Let $\A:\D(\A)\to Z$ be the generator of 
a strongly continuous group of bounded operators on $Z$.
Let $Z_{1}$ denote $\D(\A)$ with the norm $\norm{z}_{1}=\norm{(\A-\beta)z}$
for some $\beta\notin \spec(\A)$
($\spec(\A)$ denotes the spectrum of $\A$, 
this norm is equivalent to the graph norm 
and $Z_{1}$ is densely and continuously embedded in $Z$)
and let $Z_{-1}$ be the completion of $Z$ with respect to the norm 
$\norm{\zeta}_{-1}=\norm{(\A-\beta)^{-1}\zeta}$.
Let $Z'$ denote the dual of $Z$ 
with respect to the pairing $\scd{\cdot}{\cdot}$.
The dual of $\A$ is a self-adjoint operator $\A'$ on $Z'$.
The dual of $Z_{1}$ is the space $Z'_{-1}$ 
which is the completion of $Z'$ with respect to the norm 
$\norm{\zeta}_{-1}=\norm{(\A'-\bar{\beta})^{-1}\zeta}$
and the dual of $Z_{-1}$ is the space $Z'_{1}$
which is $\D(\A')$ 
with the norm $\norm{z}_{1}=\norm{(\A'-\bar{\beta}z}$.

Let $\CC\in \L(Z_{1},\V)$ and let $\CC'\in\L(\V',Z'_{-1})$ denote its dual.
Note that the same theory applies to any $\A$-bounded operator $\CC$ 
with a domain invariant by $(e^{t\A})_{t\geq 0}$ 
since it can be represented by an operator in $\L(Z_{1},\V)$ (cf.~\cite{Wei89}).

We consider the dual observation and control systems
with output function $v$ and input function $u$:
\begin{gather}
\dz(t)=\A z(t), \quad z(0)=z_{0}\in Z, 
\quad v(t)=\CC z(t), \label{eqsystz}\\
\dzeta(t)=\A' \zeta(t)+\CC' u(t), \quad \zeta(0)=\zeta_{0}\in Z', 
\quad u\in L^{2}\loc(\R;Z') \label{eqsystzeta}
. \end{gather}
  
We make the following equivalent admissibility assumptions 
on the observation operator $\CC$ and the control operator $\CC'$
(cf.~\cite{Wei89}): $\forall T>0$, $\exists K_{T}>0$,
\begin{gather}
\label{eqadmobs}
\forall z_{0}\in \D(\A),
\quad
\int_{0}^{T}\norm{\CC e^{t\A}z_{0}}^{2}dt \leq K_{T}\norm{z_{0}}^{2} , \\
\label{eqadmcon}
\forall u\in L^{2}(\R;\V'), \quad
\norm{\int_{0}^{T}e^{t\A'}\CC' u(t)dt}^{2}
\leq K_{T} \int_{0}^{T}\norm{u(t)}^{2}dt .
\end{gather} 
With this assumption, the output map 
$z_{0} \mapsto v$ 
from $\D(\A)$ to $L^{2}\loc(\R;\V)$ 
has a continuous extension to $Z$. 
The equations (\ref{eqsystz}) and (\ref{eqsystzeta})
have unique solutions $z\in C(\R,Z)$ and $\zeta\in C(\R,Z')$ 
defined by:
\begin{gather}
\label{eqmildsol}
z(t)=e^{t\A}z_{0},\quad 
\zeta(t)=e^{t\A'}\zeta(0)+\int_{0}^{t}e^{(t-s)\A}\B u(s) ds .
\end{gather}

The following dual notions of observability and controllability 
are equivalent (cf.~\cite{DR77}).
\begin{definition}
\label{def:obs} 
The system (\ref{eqsystz}) 
is {\em final observable} in time $T>0$ at cost $\kappa_{T}>0$
if the following observation inequality holds:
$\dis \forall z_{0}\in Z,
\quad
\norm{z(T)}^{2}\leq \kappa_{T}^{2}
\int_{0}^{T}\norm{v(t)}^{2}dt  $.
The system (\ref{eqsystzeta}) 
is {\em null-controllable} in time $T>0$ at cost $\kappa_{T}>0$ if 
for all $\zeta_{0}$ in $Z'$, 
there is a $u$ in $L^{2}(\R;\V')$ such that 
$\zeta(T)=0$ and
$\dis \int_{0}^{T}\norm{u(t)}^{2}dt\leq \kappa_{T}^{2}\norm{\zeta_{0}}^{2}$.
The {\em null-controllability cost} for (\ref{eqsystzeta}) 
in time $T$
is the smallest constant in the latter inequality 
(equivalently in the former observation inequality), 
still denoted $\kappa_{T}$.
When (\ref{eqsystzeta}) is not null-controllable in time $T$, 
we set $\kappa_{T}=+\infty$.
\end{definition}


\subsection{Tensor products}
\label{sec:tensprod}

Now, we introduce the specific tensor product structure 
of the abstract control systems (\ref{eqsystzeta}) 
under consideration here.
Let $X$, $Y$, $V$ be separable Hilbert spaces
and $I$ denote the identity operator on each of them. 
Let $A:D(A)\to X$ and $B:D(B)\to Y$ be generators of
strongly continuous semigroups of bounded operators on $X$ and $Y$.
Let $C\in\L(X_{1},V)$ be admissible for the control system:
\begin{gather}
\dxi(t)=A' \xi(t)+C' u(t), \quad \xi(0)=\xi_{0}\in X', 
\quad u\in L^{2}\loc(\R;V') \label{eqsystxi}
. \end{gather}
Let $X\ctensprod Y$ and $V\ctensprod Y$ denote the closure of 
the algebraic tensor products $X\otensprod Y$ and $V\otensprod Y$ 
for the natural Hilbert norms.
The operator $C\otensprod I:D(C)\otensprod Y \to V\ctensprod Y$ 
is densely defined on $X\ctensprod Y$.
The operator $A\otensprod I+I\otensprod B$ defined on 
the algebraic $D(A)\otensprod D(B)$ is closable
and its closure, denoted $A+B$, 
generates a strongly continuous semigroup of bounded operators on 
$X\ctensprod Y$.

\begin{lemma}
\label{lem:abstprod}
Let $Z=X\ctensprod Y$, $\V=V\ctensprod Y$, $\A=A+B$ and $\CC=C\otensprod I$.
If $B$ is a non-positive self-adjoint operator, 
then, for all $T> 0$, 
the null-controllability cost $\kappa_{T}$ 
for (\ref{eqsystzeta})
is lower then the null-controllability cost $k_{T}$  
for (\ref{eqsystxi}) in the same time $T$.
\end{lemma}

\begin{proof}
We may assume that $k_{T}$ is finite. By definition it satisfies:
\begin{gather} 
\label{eqobsx}
\forall x\in X,
\quad
\norm{e^{TA}}^{2}\leq k_{T}^{2}
\int_{0}^{T}\norm{C e^{tA}}^{2}dt 
. \end{gather}
We have to prove that: 
\begin{gather} 
\label{eqobsprod}
\forall z\in X\ctensprod Y
\quad
\E:=\norm{e^{T(A+B)}z}^{2}\leq k_{T}^{2}
\int_{0}^{T} \norm{(C\otensprod I)e^{t(A+B)}z}^{2} dt =:\O 
. \end{gather}

As explained in the proof of lem.~7.1 in \cite{LMconscost}:
\begin{gather} 
\label{eqsemiprod}
\forall t\geq 0, \quad e^{t(A+B)}=e^{tA}\otensprod  e^{tB} 
\ . \end{gather}                                

Applying the spectral theorem for unbounded self-adjoint operators 
on separable Hilbert spaces to $B\leq 0$ 
(cf. th.~VIII.4 in \cite{RS}),
yields a measure space $(M,\mathcal{M},\mu)$
with finite measure $\mu$,
a measurable function $b:M\to (-\infty,0]$ 
and a unitary operator $U:Y\to L^{2}(M,d\mu)$
such that:
\begin{gather} 
\label{eqU}
\forall y\in Y, \quad 
\norm{e^{tB}y}^{2}=\int_{M} e^{2tb(m)}  \abs{Uy(m)}^{2} \mu(dm)
\ .\end{gather}
Since $X$ is separable, there is a unique isomorphism 
from $X\ctensprod L^{2}(M,d\mu)$ to $L^{2}(M,d\mu; X)$
so that $x\otensprod f(m) \mapsto f(m)x$
(cf. th.~II.10 in~\cite{RS}).
We denote by $\U:X\ctensprod Y \to L^{2}(M,d\mu; X)$ 
the composition of this isomorphism with $I\otensprod U$.
Similarly, 
there is a unique isomorphism 
from $V\ctensprod L^{2}(M,d\mu)$ to $L^{2}(M,d\mu; V)$
so that $v\otensprod f(m) \mapsto f(m)v$.
We denote by $\V:V\ctensprod Y \to L^{2}(M,d\mu; V)$ 
the composition of this isomorphism with $I\otensprod U$.
By decomposing into an orthonormal basis of $X$, (\ref{eqU}) implies:  
\begin{gather} 
\label{eqUU}
\forall z\in X\ctensprod Y , \quad 
\norm{(I\otensprod e^{tB})z}^{2}=\int_{M} e^{2tb(m)}  \abs{\U z(m)}^{2} \mu(dm)
\\ \label{eqV}
\forall w\in V\ctensprod Y , \quad 
\norm{(I\otensprod e^{tB})w}^{2}=\int_{M} e^{2tb(m)}  \abs{\V w(m)}^{2} \mu(dm)
\ .\end{gather}

Let $z \in X\ctensprod Y$. 
Applying (\ref{eqobsx}) to $\U z(m)$ for fixed $m\in M$
and integrating yields:
\begin{gather*} 
\int_{M} \norm{e^{TA}\U z(m)}^{2} e^{2tb(m)} \mu(dm)
\leq k_{T}^{2}
\int_{M} e^{2Tb(m)} \int_{0}^{T} \norm{C e^{tA}\U z(m)}^{2} dt\, \mu(dm) 
\ . \end{gather*}
Since $e^{TA}\U z=\U (e^{TA}\otensprod I)z$, 
(\ref{eqUU}) and (\ref{eqsemiprod}) imply that 
the left hand side is $\E$ defined in (\ref{eqobsprod}).
Using 
Fubini's theorem and $b\leq 0$ to bound the right hand side from above yields:
\begin{gather*} 
\E \leq k_{T}^{2}
\int_{0}^{T} \int_{M}  e^{2tb(m)} \norm{C e^{tA}\U z(m)}^{2} \mu(dm)\, dt  
\ . \end{gather*}
Since $C e^{tA}\U z=\V (C e^{tA}\otensprod I)z$,
(\ref{eqV}) and (\ref{eqsemiprod}) imply  that 
the right hand side is $\O$ defined in (\ref{eqobsprod}),
which completes the proof of (\ref{eqobsprod}).
\end{proof}

\subsection{Proof of th.\ref{th:strip}, 
th.\ref{th:prodgamma} and th.\ref{th:prodomega}}
\label{sec:appli}

The first part of th.\ref{th:strip} 
is a particular case of th.\ref{th:prodgamma}.
The second part is an estimate on the null-controllability cost 
which results from def.\ref{def:alpha} and lem.\ref{lem:abstprod}
with $X=L^{2}(0,L)$, $Y=L^{2}(\R)$, $Z=\R$, 
$A=\partial_{x}^{2}$, $\D(A)=H^{2}(0,L)\cap H^{1}_{0}(0,L)$,
$B=\partial_{y}^{2}$, $\D(B)=H^{2}(\R)$,
$Cf=\partial_{x}f\res{x=0}$.
The reader balking at the abstraction of lem.\ref{lem:abstprod}
can prove it in this particular case 
using the Fourier transform on the real line in the $y$ variable 
where the spectral theorem was used
(then $\mu$ is the Lebesgue measure and $b(m)=-\abs{m}^{2}$)
and a discrete Fourier decomposition on the interval in the $x$ variable. 

Th.\ref{th:prodgamma} and th.\ref{th:prodomega}
are direct applications of lem.\ref{lem:abstprod}
with $X=L^{2}(M)$, $Y=L^{2}(\Mt)$, 
$A=\Delta$, $\D(A)=H^{2}(M)\cap H^{1}_{0}(M)$,
$B=\Deltat$, $\D(B)=H^{2}(\Mt)\cap H^{1}_{0}(\Mt)$.
Th.\ref{th:prodgamma} corresponds to 
$Z=L^{2}(\Gamma)$ and $Cf=\partial_{\nu}f\res{\Gamma}$
where $\partial_{\nu}$ denotes 
the exterior Neumann vector field on $\partial M$.
Th.\ref{th:prodomega} corresponds to 
$Z=L^{2}(\Omega)$ and $Cf=f\res{\Omega}$.


\section{Geometric necessary condition.}
\label{sec:GNC}

In this section, we prove th.\ref{th:GNC}.
Henceforth, the domain of the Laplacian is 
$\D(\Delta)=H^{2}(M)\cap H^{1}_{0}(M)$.
Since controllability and observability in def.\ref{def:obs} are equivalent, 
the heat equation (\ref{eqHeatI}) is null-controllable in time $T$ 
if and only if there is a $C_{\Omega,T}>0$ such that 
\begin{gather}
\label{eqobs}
\forall f_{0}\in L^{2}(M),\quad 
\int_{M}\abs{e^{T\Delta}f_{0}}^{2}dx\leq  C_{\Omega,T}
\int_{0}^{T}\int_{M} \abs{e^{t\Delta}f_{0}}^{2}dx\, dt
\ .
\end{gather}

As for th.2.1 in \cite{LMheatcost}
where the null-controllability cost $C_{\Omega,T}$ (on a compact $M$) 
was bounded from below as $T\to 0$,
the strategy is to choose the initial datum $f_{0}$ 
to be an approximation of the Dirac mass $\delta_{y}$ at some $y\in M$ 
which is as far from $\Omega$ as possible.
Therefore both proofs build on heat kernel estimates.
But here we need estimates which are uniform on $M$ for compact times 
and we use the finer notion of averaged distance of $y$ to $\Omega$ 
(cf. def.\ref{def:dist}).

\subsection{Heat kernel estimates.}
Let $K_{M}(t,x,y)$ denote the Dirichlet heat kernel on $M$
(i.e.~the fundamental solution ``$e^{t\Delta}\delta_{y}(x)$'').
We recall some well-known facts about it.
The heat kernel on $M$ satisfies the following upper bound
(cf.~th.3.2.7 in \cite{Dav89}): 
$\forall \eps\in ]0,1[$, $\exists a_{\eps}>0$ s.t. 
\begin{gather}\label{equb}
\forall t>0,\, \forall x,y\in M, \quad K_{M}(t,x,y) \leq 
a_{\eps}t^{-n/2}\exp\left(-\frac{d(x,y)^{2}}{4(1+\eps)t}\right) \ .
\end{gather}
Let $C$ be a bounded open subset of $M$. 
Let $(\lambda_{j})_{j\in \mathbb{N}^{*}}$
be a nondecreasing sequence of nonnegative real numbers 
and $(e_{j})_{j\in \mathbb{N}^{*}}$
be an orthonormal basis of $L^{2}(M)$
such that $e_{j}$ is an eigenfunction of the Dirichlet Laplacian on $C$
with eigenvalue $-\lambda_{j}$.
By the maximum principle, 
the heat kernel on $M$ satisfies the lower bound:
\begin{gather}\label{eqlb}
\forall t>0,\, \forall x,y\in C, \quad K_{M}(t,x,y)\geq K_{C}(t,x,y) 
=\sum_{j}e^{-t\lambda_{j}}
e_{j}(y)e_{j}(x) \ .
\end{gather}

From these pointwise bounds on the heat kernel,
we deduce bounds for the $L^{2}$ norms appearing in (\ref{eqobs}).
Def.\ref{def:dist} and (\ref{equb}) imply
\begin{gather}\label{equbL}
\int_{T_{1}}^{T_{2}}
\int_{\Omega}\abs{K_{M}(t,x,y))}^{2} dx\, dt
\leq 
a_{\eps}^{2} \frac{T_{2}-T_{1}}{T_{1}^{n}}
\exp\left(-\frac{\dbar_{(1+\eps)T_{2}}(y,\Omega)^{2}}{2(1+\eps)T_{2}}\right)
\ .
\end{gather}
If $C\subset M$ is an $n$-dimensional cube 
with center $y$ and half diagonal length $d$,
i.e. with edge length $c=2d/\sqrt{n}$,
then the first eigenvalue and eigenfunction 
of the Dirichlet Laplacian on $C$ are
$\dis \lambda_{1}=n\left(\frac{\pi}{2c}\right)^{2}$ and  
$\dis e_{1}(x)
=c^{-n/2}\prod_{m=1}^{n}\cos\left(\frac{\pi (x_{m}-y_{m})}{2c}\right)$.
Therefore, (\ref{eqlb}) imply
\begin{gather}\label{eqlbL}
\int_{M}\abs{K_{M}(t,x,y))}^{2} dx 
\geq \int_{C}\abs{K_{C}(t,x,y))}^{2} dx 
\geq e^{-2\lambda_{1} t}\abs{e_{1}(y)}^{2}
= \frac{n^{n/2}}{(2d)^{n}} \exp\left(-\frac{\pi^{2}n^{2}t}{8d^{2}}\right).
\end{gather}

\begin{remark}
We tried without tangible improvement 
to deduce $L^{2}$ lower bounds on the heat kernel 
from the uniform pointwise lower bounds 
available in the literature (cf.~\cite{VDB90}) 
instead of deducing it from the more basic fact (\ref{eqlb}).
\end{remark}

\subsection{Proof of th.\ref{th:GNC}} 

Let $\set{y_{k}}_{k\in\N}$, $\Tbar$ and $\kappa$ 
satisfy the geometric condition (\ref{eqgeomcond}).
By contradiction, assume that 
the heat equation (\ref{eqHeatI}) is null-controllable in some time $T<\Tbar$,
i.e. the observability inequality (\ref{eqobs}) holds for some $C_{\Omega,T}$.
Let $\eps\in ]0,1[$, $\eps<\kappa-1$, and let $\kappa'=\kappa (1+\eps)^{-1}>1$.
Let $\alpha>0$ be such that $\Tbar=(1+\alpha)(1+\eps)T$
and let $\Tdbar=(1+\alpha)T$.
Since $\ddbar_{T}/T$ is non-increasing, 
(\ref{eqgeomcond}) implies
\begin{gather}
\label{eqgc}
s_{k}:=\frac{\dbar_{\Tbar}(y_{k}, \Omega)^{2}}{2\Tbar}
-\kappa'\frac{\pi^{2}n^{2}\Tdbar}{8\ddbar_{\Tdbar}(y_{k},\partial M)^{2}}
\to +\infty\, , \text{ as } k\to +\infty 
\ .
\end{gather}
Let $k\in\N$ and let $f_{0}(x)=K_{M}(\alpha T, x,y_{k})$
so that $e^{t\Delta}f_{0}(x)=K_{M}(\alpha T + t, x,y_{k})$.
Plugging into (\ref{eqobs}) the upper bound 
(\ref{equbL}) with $T_{1}=\alpha T$ and $T_{2}=\Tdbar$
and the lower bound (\ref{eqlbL}) for the cube $C$ with center $y_{k}$
and half diagonal length $d=\ddbar_{\Tdbar}(y_{k},\partial M)$
(this is just the optimal choice for $d$)
yields:
\begin{gather*}
\frac{n^{n/2}}{(2\ddbar_{\Tdbar}(y_{k},\partial M))^{2}}
\exp\left(-\frac{\pi^{2}n^{2}\Tdbar}{8\ddbar_{\Tdbar}(y_{k},\partial M)^{2}}\right)
\leq C_{\Omega,T}   \frac{a_{\eps}^{2}}{\alpha^{n}T^{n-1}}
\exp\left(-\frac{\dbar_{\Tbar}(y_{k},\Omega)^{2}}{2\Tbar}\right)
\ .
\end{gather*}
Since $\kappa'>1$, we deduce that there is an $s>0$ independent of $k$ 
such that $\ln C_{\Omega,T}\geq s_{k}- s$ 
and $\lim_{k}s_{k}=+\infty$ as in (\ref{eqgc}).
This contradicts the existence of $C_{\Omega,T}$ 
and completes the proof of th.\ref{th:GNC}.

\subsection{Proof of th.\ref{th:rod} iii) and another example} 
\label{sec:proof}

To prove that the geometric condition (\ref{eqgeomcond}) holds 
for $M$ and $\Omega$ defined in th.\ref{th:rod} iii),
we consider a sequence $m_{k}=(0,0,z_{k})\in M$ with $\lim_{k}z_{k}=+\infty$.
Since $S$ is bounded, we may assume that $R$ is bounded.
Let $G_{T}(z)=\exp(-z^{2}/(2T))$ and  
let $D(z)$ denote the disk with center $(0,0)$ and radius $R(z)$.
We have: 
\begin{gather*}
\begin{split}
  I_{k}&:=\int_{\Omega}G_{T}(d(m_{k},m)) dm
=\int_{\R}\exp\left(-\frac{(z-z_{k})^{2}}{2T}\right)
\int_{D(z)}\exp\left(-\frac{x^{2}+y^{2}}{2T}\right) dx\, dy\, dz \\ &
\leq \int_{\R}\pi R(z)^{2}G(z-z_{k})dz
= \pi R^{2}\ast G_{T} (z_{k}) \to 0\, , \text{ as } k\to +\infty \,  ,
\end{split}
\end{gather*}
since $G_{T}\in L^{1}(\R)$, $R^{2}\in L^{\infty}(\R)$ and 
$\lim_{\abs{z}\to \infty}R(z)=0$. 
Therefore, by def.\ref{def:dist}, 
$\dbar_{\Tbar}(m_{k}, \Omega)^{2}=-2T\ln I_{k}\to +\infty$
and, since $(0,0)\in S$, 
$\ddbar_{\Tbar}(m_{k},\partial M)^{2}\geq d_{\partial}(m_{k})^{2}
=\inf_{(x,y)\in \R^{2}\setminus S}(x^{2}+y^{2})>0$.
Hence (\ref{eqgeomcond}) holds for the sequence
$\set{m_{k}}_{k\in\N}$ with any $\Tbar$ and $\kappa$,
which completes the proof of th.\ref{th:rod} iii).

\begin{remark}\label{rem:shrinkrod}
To illustrate the usefulness of 
the second term in the geometric condition (\ref{eqgeomcond}),
we give an example close to th.\ref{th:rod} ii)
where (\ref{eqgeomcond}) is satisfied by 
a sequence $\set{m_{k}}_{k\in\N}$ tending to the boundary of $M$.

Consider the shrinking rod 
$M=\set{(x,y,z)\in\R^{3}\st x^{2}+y^{2}< R(\abs{z})^{2}}$
where the continuous non-increasing 
function $R:[0,\infty)\to ]0,\infty)$ tends to zero at infinity.
The heat equation (\ref{eqHeatI}) is not null-controllable in any time $T>0$
from any interior region $\Omega$ of finite Lebesgue measure
such that 
$M\setminus\Omega$ contains  a sequence of slabs 
$S_{k}:=\set{(x,y,z)\in\R^{2}\times[0,\infty)\st 
x^{2}+y^{2}< R(z)^{2}, \abs{z-z_{k}}\leq d_{k} }$ 
satisfying 
\begin{gather*}
\exists\kappa'>1,\quad 
d_{k}^{2}-
\kappa' 
\frac{\pi^{2}n^{2}}{4} 
\left(\frac{T}{R(z_{k}+d_{k})}\right)^{2}
\to +\infty\, , \text{ as } k\to +\infty 
\ .
\end{gather*}
Indeed $m_{k}=(0,0,z_{k})$ satisfies 
$d_{\partial}(m_{k})\geq R(z_{k}+d_{k})$ for $d_{k}\geq \norm{R}_{L^{\infty}}$,
and $d(m_{k},\Omega)\geq d_{k}$.
Hence  $\set{m_{k}}$ satisfies (\ref{eqgeomcond})
for any $\kappa\in]1,\kappa'[$ and $\Tbar=\sqrt{\kappa'/\kappa}T>T$.
In particular, if 
$\lim_{z\to +\infty}z R(z)=+\infty$ (i.e. $M$ does not shrink too fast)
then the heat equation (\ref{eqHeatI}) is not null-controllable in any time $T$
from any bounded 
$\Omega$.
\end{remark}



\def\cprime{$'$}
\providecommand{\bysame}{\leavevmode\hbox to3em{\hrulefill}\thinspace}
\providecommand{\MR}{\relax\ifhmode\unskip\space\fi MR }
\providecommand{\MRhref}[2]{%
  \href{http://www.ams.org/mathscinet-getitem?mr=#1}{#2}
}
\providecommand{\href}[2]{#2}

\end{document}